\documentclass[reqno]{amsart}
\usepackage{graphicx}
\usepackage{fullpage}

\theoremstyle{definition}

\newcommand{\taumin}{\tau_{\rm min}}
\newcommand{\taumax}{\tau_{\rm max}}
\newcommand{\xmin}{x_{\rm min}}
\newcommand{\xmax}{x_{\rm max}}
\newcommand{\green}{f_{_{\rm G}}}
\newcommand{\tauperp}{\tau_\perp}

\newcommand{\deltapar}{\beta}
\newcommand{\chimin}{\chi_{\rm min}}
\newcommand{\chimax}{\chi_{\rm max}}

\def\ngreen{n_{_{\rm G}}}
\def\ugreen{U_{_{\rm G}}}
\def\sigperp{\sigma_\perp}
\def\sigpar{\sigma_\|}
\def\sigbar{\overline\sigma}

\def\colrad{r_0}

\theoremstyle{remark}

\begin{document}

\title[Spectral formation in x-ray pulsars]
{Spectral formation in x-ray pulsars and associated identities involving
the Laguerre polynomials}

\author[Peter A. Becker]{Peter A. Becker\\
Department of Computational and Data Sciences,\\
College of Science,\\
George Mason University,\\
Fairfax, VA 22030-4444, USA}

\email{pbecker@gmu.edu}

\subjclass{Primary 33C10, 33C45, 34B27; Secondary 85A25}

\date{Submitted June 2, 2007.}

\keywords{Confluent Hypergeometric Functions, Green's Functions,
Orthogonal Polynomials, Radiative Transfer}

\vskip2.0truein
\centerline{accepted for publication in the Journal of Mathematical Physics}

\begin{abstract}
The partial differential equation governing the formation of the
radiation spectrum in the shocked plasma accreting onto an x-ray pulsar
is solved to obtain the associated Green's function describing the
scattering of monochromatic radiation injected into the plasma at a
fixed altitude above the star. Collisions between the photons and the
infalling electrons cause both the ordered and random components of the
plasma energy to be transferred to the radiation, which escapes by
diffusing through the walls of the accretion column. The analytical
solution for the Green's function provides important physical insight
into the formation of the observed spectra in x-ray pulsars. Interesting
mathematical aspects of the problem include the establishment of a new
summation formula involving the Laguerre polynomials, based on the
calculation of the photon number density via integration of the Green's
function.
\end{abstract}

\maketitle

\bigskip

\section*{\bf I. INTRODUCTION}

\bigskip

The analysis of solutions to partial differential equations satisfying
physically prescribed boundary conditions often leads to new insights
into the relationships between the various special functions of
mathematical physics. In this article, we employ the methods of
classical analysis to obtain the closed-form solution for the Green's
function describing the scattering of radiation inside the plasma
accreting onto the surface of a rotating neutron star in an x-ray
pulsar. The method is based on the separation of the two-dimensional
partial differential equation into two second-order ordinary
differential equations depending on the energy and space variables. By
applying suitable boundary conditions in the spatial dimension, we are
able to determine the eigenvalues and show that the corresponding
spatial eigenfunctions are Laguerre polynomials. The solutions to the
associated energy equation are Whittaker functions. The global solution
for the two-dimensional Green's function is obtained by forming an
expansion based on the energy and space eigenfunctions. The availability
of this solution also allows us to calculate the distribution of the
radiation escaping from the accreting gas, which forms the spectrum
observed at earth. Beyond the direct physical relevance of the Green's
function, the method of solution also yields several additional results
of mathematical interest, including a new summation identity for the
Laguerre polynomials.

We provide a brief overview of the physical problem before proceeding
with the main derivation. The fundamental power source for the radiation
produced in bright x-ray pulsars is the gravitational accretion (inflow)
of fully ionized plasma that is channeled onto the poles of the rotating
neutron star by the strong magnetic field. In these luminous sources,
the pressure of the photons governs the dynamical structure of the
accretion flow, and therefore the gas must pass through a
radiation-dominated shock on its way to the stellar surface. The kinetic
energy of the gas is carried away by the high-energy radiation that
escapes from the column, which allows the material to come to rest on
the stellar surface.$^{1,2}$ The strong compression of the infalling gas
drives its temperatures up to several million Kelvins. The gas therefore
radiates x-rays, which appear to pulsate due to the star's rotation. The
observed x-ray spectrum is often distinctly nonthermal, indicating that
nonequilibrium processes are playing an important role in the radiative
transfer occurring inside the accretion column.$^{3}$

The nonthermal shape of the spectrum is due to the transfer of energy
from the gas to the photons via electron scattering. The plasma
possesses both ``ordered'' kinetic energy associated with the inflow,
and ``random'' kinetic energy associated with the thermal motion of the
particles. These two types of energy are transferred to the photons via
first- and second-order Fermi processes, respectively. Our primary goal
in this article is to obtain the exact solution for the Green's function
describing the effect of electron scattering on monoenergetic seed
photons injected from a source located at an arbitrary altitude inside
the accretion column. The Green's function contains a complete
representation of the fundamental physics governing the propagation of
the photons in the physical and energy spaces. Since the transport
equation governing the radiation spectrum is linear, we can compute the
solution associated with any seed photon distribution via convolution.
In x-ray pulsars, the most important sources of seed photons are
bremsstrahlung, cyclotron, and blackbody emission.$^{4}$ The
Green's function provides a direct means for investigating how photons
produced by the various source mechanisms contribute to the observed
nonthermal spectra in x-ray pulsars.

\bigskip

\section*{\bf II. FUNDAMENTAL EQUATIONS}

\bigskip

The plasma flowing onto the neutron star is assumed to have a
cylindrical geometry, maintained by the strong dipole magnetic field. We
define the spatial coordinate $z$ as the altitude measured from the
stellar surface along the axis of the cylindrical accretion column. The
gas flows onto the star with velocity $v(z) < 0$, which vanishes at the
surface of the star. The Green's function, $\green(z_0,z,\epsilon_0,
\epsilon)$, is defined as the radiation distribution measured at
location $z$ and energy $\epsilon$ resulting from the injection of $\dot
N_0$ photons per second with energy $\epsilon_0$ from a monochromatic
source located at altitude $z_0$ inside the column. In a steady-state
situation, $\green$ satisfies the Kompaneets transport equation$^{4,5}$
\begin{eqnarray}
v \, {\partial \green \over \partial z}
&=& {dv \over d z}\,{\epsilon\over 3} \,
{\partial \green\over\partial\epsilon}
+ {\partial\over\partial z}
\left({c\over 3 n_e \sigpar}\,{\partial \green\over\partial z}\right)
- {\green \over t_{\rm esc}}
\nonumber
\\
&+& {n_e \sigbar c \over m_e c^2} {1 \over\epsilon^2}
{\partial\over\partial\epsilon}\left[\epsilon^4\left(\green
+ k T_e \, {\partial \green\over\partial\epsilon}\right)\right]
+ {\dot N_0 \, \delta(\epsilon-\epsilon_0) \, \delta(z-z_0)
\over \pi r_0^2 \epsilon_0^2}
\ ,
\label{eq2.1}
\end{eqnarray}
where $r_0$ is the radius of the accretion column, $t_{\rm esc}$
represents the mean time photons spend in the plasma before escaping
through the walls of the column, $\sigbar$ denotes the angle-averaged
electron scattering cross section, $c$ is the speed of light, $k$ is
Boltzmann's constant, and $T_e$, $n_e$, and $m_e$ denote the electron
temperature, number density, and mass, respectively. The left-hand side
of (\ref{eq2.1}) represents the comoving time derivative of the Green's
function, and the terms on the right-hand side represent first-order
Fermi energization (``bulk Comptonization''), spatial diffusion along
the column axis, photon escape, stochastic (thermal) Comptonization, and
photon injection, respectively.

The transport equation employed here is similar to the one analyzed by
Becker$^{1}$ and Becker and Wolff,$^{2}$ except that thermal
Comptonization has now been included via the appearance of the
Kompaneets (1957) operator, which is the fourth term on the right-hand
side of (\ref{eq2.1}). The radiation number density $\ngreen$ and energy
density $\ugreen$ associated with the Green's function are given by
\begin{equation}
\ngreen(z) = \int_0^\infty \epsilon^2 \,
\green(z_0,z,\epsilon_0,\epsilon) \, d\epsilon \ , \ \ \ \ \ 
\ugreen(z) = \int_0^\infty \epsilon^3 \,
\green(z_0,z,\epsilon_0,\epsilon) \, d\epsilon
\ .
\label{eq2.2}
\end{equation}
The Green's function provides fundamental physical insight into the
spectral redistribution process, and it also allows us to calculate the
particular solution for the spectrum $f(z,\epsilon)$ associated with an
arbitrary photon source $Q(z,\epsilon)$ using the integral
convolution$^{6}$
\begin{equation}
f(z,\epsilon) = \int_0^\infty\int_0^\infty
{\green(z_0,z,\epsilon_0,\epsilon) \over \dot N_0} \ \epsilon_0^2
\, Q(z_0,\epsilon_0) \, d\epsilon_0 \, dz_0
\ ,
\label{eq2.3}
\end{equation}
where the source function $Q$ is normalized so that $\epsilon^2 \,
Q(z,\epsilon) \, d\epsilon \, dz$ gives the number of seed photons
injected per unit time in the altitude range between $z$ and $z+dz$ and
the energy range between $\epsilon$ and $\epsilon+d\epsilon$.

Following Lyubarskii and Sunyaev,$^{7}$ we will assume that the
electron temperature $T_e$ has a constant value, which is physically
reasonable since most of the Comptonization occurs in a relatively
compact region near the base of the accretion column. In this situation,
it is convenient to work in terms of the dimensionless energy variable
$\chi$, defined by
\begin{equation}
\chi(\epsilon) \equiv {\epsilon \over kT_e}
\ ,
\label{eq2.4}
\end{equation}
and the dimensionless optical depth variable $\tau$, defined by
\begin{equation}
d\tau \equiv n_e(z) \, \sigpar \, dz
\ , \ \ \ \
\tau(z) \equiv \int_0^z n_e(z') \, \sigpar \, dz'
\ ,
\label{eq2.5}
\end{equation}
where $z$ and $\tau$ both vanish at the stellar surface. The electron
number density $n_e$ is related to the inflow velocity $v$ via
\begin{equation}
\dot M \equiv \pi r_0^2 \, m_p \, n_e \, |v| = {\rm constant}
\ ,
\label{eq2.6}
\end{equation}
where $\dot M$ is the mass accretion rate onto the magnetic pole of the
star and $m_p$ is the proton mass. Note that we can write the Green's
function as either $\green(z_0,z, \epsilon_0,\epsilon)$ or
$\green(\tau_0,\tau,\chi_0,\chi)$ since the variables $(z,\epsilon)$ and
$(\chi,\tau)$ are interchangeable by utilizing (\ref{eq2.4}) and
(\ref{eq2.5}).

Making the change of variable from $z$ to $\tau$ in (\ref{eq2.1}), we
find after some algebra that the transport equation for the Green's
function can be written in the form
\begin{eqnarray}
{v \over c} \, {\partial \green \over \partial \tau}
&=& {1 \over c}\,{dv \over d\tau}\,{\chi\over 3} \,
{\partial \green\over\partial\chi}
+ {1 \over 3}\,{\partial^2\green\over\partial\tau^2}
- {\xi^2\,v^2 \over c^2} \, \green
\nonumber
\\
&+& {\sigbar \over \sigpar} \, {kT_e \over m_e c^2}
{1 \over\chi^2}{\partial\over\partial\chi}\left[\chi^4\left(\green
+ {\partial\green\over\partial\chi}\right)\right]
+ {\dot N_0 \, \delta(\chi-\chi_0) \, \delta(\tau-\tau_0)
\over \pi r_0^2 c k T_e \epsilon_0^2}
\ ,
\label{eq2.7}
\end{eqnarray}
where $\tau_0 \equiv \tau(z_0)$, $\chi_0 \equiv \chi(\epsilon_0)$, and
we have introduced the dimensionless constant
\begin{equation}
\xi \equiv {\pi r_0 m_p c \over \dot M (\sigpar \sigperp)^{1/2}}
\ ,
\label{eq2.8}
\end{equation}
which determines the importance of the escape of photons from the
accretion column. It can be shown that $\xi$ is roughly equal to the
ratio of the dynamical (accretion) timescale divided by the timescale
for the photons to diffuse through the column walls,$^{4}$ and therefore
the condition $\xi \sim 1$ must be satisfied in order to ensure that the
radiation pressure decelerates the gas to rest at the stellar surface,
with the kinetic energy of the material carried away by the photons
escaping through the column walls.$^{8}$

Following Becker and Wolff,$^{2}$ we compute the mean escape time
using the diffusive prescription
\begin{equation}
t_{\rm esc}(z) = {r_0 \, \tauperp \over c} \ , \ \ \ \ \
\tauperp(z) = n_e \, \sigperp \, \colrad
\ ,
\label{eq2.9}
\end{equation}
where $\tauperp$ represents the perpendicular scattering optical
thickness of the cylindrical accretion column. Note that $\tauperp$ and
$t_{\rm esc}$ are each functions of $z$ through their dependence on the
electron number density $n_e(z)$. Becker$^{8}$ confirmed that the
diffusion approximation employed in (\ref{eq2.9}) is valid since
$\tauperp > 1$ for typical x-ray pulsar parameters. The technical
approach used to solve for the Green's function, carried out in \S~3,
involves the derivation of eigenvalues and associated eigenfunctions
based on the set of spatial boundary conditions for the problem.$^{9,
10,11,12}$

\bigskip

\section*{\bf III. EXACT SOLUTION FOR THE GREEN'S FUNCTION}

\bigskip

Lyubarskii and Sunyaev$^{7}$ demonstrated that the transport
equation~(\ref{eq2.7}) is separable in energy and space if $\chi \ne
\chi_0$ and the velocity profile has the particular form
\begin{equation}
v(\tau) = - \alpha \, c \, \tau
\ ,
\label{eq3.1}
\end{equation}
where $\alpha$ is a positive constant. This expression provides a
reasonable approximation of the actual velocity profile in an x-ray
pulsar accretion column. By combining (\ref{eq2.5}), (\ref{eq2.6}),
(\ref{eq2.8}), and (\ref{eq3.1}), we can express $\tau$ as an explicit
function of the altitude $z$, obtaining
\begin{equation}
\tau(z) = \left(\sigpar \over \sigperp\right)^{1/4}
\left(2 \, z \over \alpha \, \xi \, r_0\right)^{1/2}
\ .
\label{eq3.2}
\end{equation}
Using this result to substitute for $\tau$ in (\ref{eq3.1}), we note
that the velocity profile required for separability is related to $z$
via
\begin{equation}
v(z) = - \left(\sigpar \over \sigperp\right)^{1/4}
\left(2 \, \alpha z \over \xi \, r_0\right)^{1/2} c
\ .
\label{eq3.3}
\end{equation}
This profile describes a flow that stagnates at the stellar surface
($\tau=0$, $z=0$) as required.

Adopting the velocity profile given by (\ref{eq3.1}), we find that
the transport equation~(\ref{eq2.7}) for the Green's function can be
reorganized to obtain
\begin{eqnarray}
{\alpha \, \chi\over 3} \,
{\partial \green \over \partial \chi}
- {\sigbar \over \sigpar} \, {kT_e \over m_e c^2}
{1 \over\chi^2}{\partial\over\partial\chi}\left[\chi^4\left(\green
+ {\partial\green\over\partial\chi}\right)\right]
&=& {1 \over 3}\,{\partial^2\green\over\partial\tau^2}
+ \alpha \, \tau \, {\partial \green \over \partial \tau}
- \xi^2 \alpha^2\,\tau^2 \, \green
\nonumber
\\
&+& {\dot N_0 \, \delta(\chi-\chi_0) \, \delta(\tau-\tau_0)
\over \pi r_0^2 c k T_e \epsilon_0^2}
\ .
\label{eq3.4}
\end{eqnarray}
When $\chi \ne \chi_0$, the $\delta$-function source term on the
right-hand side makes no contribution, and the differential equation is
therefore linear and homogeneous. In this case the transport equation
can be separated in energy and space using the functions
\begin{equation}
f_\lambda(\tau,\chi) = g(\lambda,\tau) \ h(\lambda,\chi) \ ,
\label{eq3.5}
\end{equation}
where $\lambda$ is the separation constant. We find that the spatial
and energy functions, $g$ and $h$, respectively, satisfy the differential
equations
\begin{equation}
{1 \over 3}\,{d^2 g \over d\tau^2}
+ \alpha \, \tau \, {d g \over d\tau}
+ \left({\alpha \lambda \over 3}
- \xi^2 \alpha^2 \tau^2\right) g = 0
\ ,
\label{eq3.6}
\end{equation}
and
\begin{equation}
{1 \over \chi^2}{d\over d\chi}\left[\chi^4
\left(h + {d h\over d\chi}\right)\right]
- \deltapar \, \chi \, {d h \over d\chi}
- \deltapar \, \lambda \, h = 0
\ ,
\label{eq3.7}
\end{equation}
where the parameter $\deltapar$ is defined by
\begin{equation}
\deltapar \equiv {\alpha \over 3} \,
{\sigpar \over \sigbar} \, {m_e c^2 \over kT_e}
\ .
\label{eq3.8}
\end{equation}
It can be shown that $\deltapar$ determines the relative importance of
bulk and thermal Comptonization. When $\deltapar$ is of order unity, the
two processes are comparable, and when $\deltapar > 1$, the bulk process
dominates.$^{4}$

\bigskip

\section*{\bf 3.1 Eigenvalues and Spatial Eigenfunctions}

\bigskip

In order to obtain the global solution for the spatial separation
function $g$, we need to understand the physical boundary conditions
that $g$ must satisfy. In the downstream region, as the gas approaches
the stellar surface, we require that the advective and diffusive
components of the radiation flux must both vanish due to the divergence
of the electron density as the gas settles onto the star. The advective
flux is indeed negligible at the stellar surface since $v \to 0$ as
$\tau \to 0$ [see Eq.~(\ref{eq3.1})]. However, in order to ensure
that the diffusive flux vanishes, we must require that $dg/d\tau \to 0$
as $\tau \to 0$. Conversely, in the upstream region, we expect that $g
\to 0$ as $\tau \to \infty$ since no photons can diffuse to large
distances in the direction opposing the plasma flow. With these boundary
conditions taken into consideration, we find that the fundamental
solution for the spatial separation function $g$ has the general form
\begin{equation}
g(\lambda,\tau) \propto
\begin{cases}
e^{-\alpha(3+w)\tau^2/4} \, U\big(\rho,{1 \over 2},{\alpha w \tau^2
\over 2}\big) \ , & \tau \ge \tau_0 \ , \cr
e^{-\alpha(3+w)\tau^2/4} M\big(\rho,{1 \over 2},{\alpha w \tau^2
\over 2}\big) \ , & \tau \le \tau_0 \ , \cr
\end{cases}
\label{eq3.9}
\end{equation}
where $M$ and $U$ denote confluent hypergeometric functions,$^{13}$
and we have made the definitions
\begin{equation}
\rho \equiv {w + 3 - 2 \lambda \over 4 w} \ ,
\ \ \ \ \
w \equiv \left(9 + 12 \, \xi^2\right)^{1/2}
\ .
\label{eq3.10}
\end{equation}

Equation~(\ref{eq3.6}) is linear, second-order, and homogeneous, and
consequently both the function $g$ and its derivative $dg/d\tau$ must be
continuous at the source location, $\tau=\tau_0$. The smooth merger of
the $M$ and $U$ functions at the source location requires that their
Wronskian, $\omega(\tau)$, must vanish at $\tau=\tau_0$, where
\begin{equation}
\omega(\tau) \equiv M\left(\rho,{1 \over 2},{\alpha w \tau^2 \over 2}\right)
{d \over d\tau}\,U\left(\rho,{1 \over 2},{\alpha w \tau^2 \over 2}\right)
- \ U\left(\rho,{1 \over 2},{\alpha w \tau^2 \over 2}\right)
{d \over d\tau}\,M\left(\rho,{1 \over 2},{\alpha w \tau^2 \over 2}\right)
\ .
\label{eq3.11}
\end{equation}
This condition can be used to solve for the eigenvalues of the
separation constant $\lambda$. By employing equation~(13.1.22) from
Abramowitz and Stegun$^{13}$ to evaluate the Wronskian, we obtain the
eigenvalue equation
\begin{eqnarray}
\omega(\tau)
= - {\Gamma({1 \over 2}) \, (2 \alpha w)^{1/2} \over \Gamma(\rho)}
\ e^{\alpha w \tau^2/2} = 0
\ .
\label{eq3.12}
\end{eqnarray}
The left-hand side vanishes when $\Gamma(\rho) \to \pm \infty$, which
implies that $\rho=-n$, where $n=0,1,2,$\ldots By combining this result
with (\ref{eq3.10}), we conclude that the eigenvalues
$\lambda_n$ are given by
\begin{equation}
\lambda_n = {4nw + w + 3 \over 2} \ ,
\ \ \ \ n=0,1,2,\ldots
\label{eq3.13}
\end{equation}

When $\lambda=\lambda_n$, the spatial separation functions $g$ reduce to
a set of global eigenfunctions $g_n$ that satisfy the boundary
conditions at large and small values of $\tau$. In this case, we can use
equations~(13.6.9) and (13.6.27) from Abramowitz and Stegun$^{13}$ to
show that the confluent hypergeometric functions $M$ and $U$ appearing
in (\ref{eq3.11}) are proportional to the generalized Laguerre
polynomials $L_n^{(-1/2)}$, and consequently the global solutions for
the spatial eigenfunctions can be written as
\begin{equation}
g_n(\tau) \equiv g(\lambda_n,\tau)
= e^{-\alpha(3+w)\tau^2/4} \, L_n^{(-1/2)}\left(\alpha w \tau^2 \over 2
\right)
\ .
\label{eq3.14}
\end{equation}
Based on equation~(7.414.3) from Gradshteyn and Ryzhik,$^{14}$ we
note that the spatial eigenfunctions satisfy the orthogonality relation
\begin{equation}
\int_0^\infty e^{3\alpha\tau^2/2} g_n(\tau) \, g_m(\tau)
\, d\tau =
\begin{cases}
{\Gamma(n+1/2) \over n! \, \sqrt{2 \alpha w}} \ , & n = m \ , \cr
0 \ , & n \ne m \ . \cr
\end{cases}
\label{eq3.15}
\end{equation}

\bigskip

\section*{\bf 3.2 Energy Eigenfunctions}

\bigskip

The solution for the energy separation function $h$ depends on the
boundary conditions imposed in the energy space. As $\chi \to 0$, we
require that $h$ not increase faster than $\chi^{-3}$ since the Green's
function must possess a finite total photon number density [see
Eq.~(\ref{eq2.2})]. Conversely, as $\chi \to \infty$, we require that
$h$ decrease more rapidly than $\chi^{-4}$ in order to ensure that the
Green's function contains a finite total photon energy density.
Furthermore, in order to avoid an infinite diffusive flux in the energy
space at $\chi=\chi_0$, the function $h$ must be continuous there. The
fundamental solution for the energy eigenfunction that satisfies the
various boundary and continuity conditions can be written as
\begin{equation}
h_n(\chi) \equiv h(\lambda_n,\chi) =
\begin{cases}
\chi^{\kappa-4} \, e^{-\chi/2} \, W_{\kappa,\mu}(\chi_0)
\, M_{\kappa,\mu}(\chi) \ ,
& \chi \le \chi_0 \ , \cr
\chi^{\kappa-4} \, e^{-\chi/2} \, M_{\kappa,\mu}(\chi_0)
\, W_{\kappa,\mu}(\chi) \ ,
& \chi \ge \chi_0 \ , \cr
\end{cases}
\label{eq3.16}
\end{equation}
where $M_{\kappa,\mu}$ and $W_{\kappa,\mu}$ denote Whittaker's functions,
and we have made the definitions
\begin{equation}
\kappa \equiv {1 \over 2} \, (\deltapar+4) \ ,
\ \ \ \ \
\mu \equiv {1 \over 2} \left[(3-\deltapar)^2 + 4 \, \deltapar \lambda_n \right]
^{1/2}
\ .
\label{eq3.17}
\end{equation}
Note that each of the eigenvalues $\lambda_n$ results in a different
value for $\mu$, and the parameter $\deltapar$ is defined in
(\ref{eq3.8}). Equation~(\ref{eq3.16}) can also be written in the more
compact form
\begin{equation}
h_n(\chi) = \chi^{\kappa-4} \, e^{-\chi/2} \, M_{\kappa,\mu}(\chimin)
\, W_{\kappa,\mu}(\chimax)
\ ,
\label{eq3.18}
\end{equation}
where
\begin{equation}
\chimin \equiv \min(\chi,\chi_0) \ ,
\ \ \ \ \
\chimax \equiv \max(\chi,\chi_0)
\ .
\label{eq3.19}
\end{equation}

\bigskip

\section*{\bf 3.3 Eigenfunction Expansion}

\bigskip

The spatial eigenfunctions $g_n(y)$ form an orthogonal set, as expected
since this is a standard Sturm-Liouville problem. The solution for the
Green's function can therefore be expressed as the infinite series
\begin{equation}
\green(\tau_0,\tau,\chi_0,\chi)
= \sum_{n=0}^\infty \ C_n \, g_n(\tau) \, h_n(\chi)
\ ,
\label{eq3.20}
\end{equation}
where the expansion coefficients $C_n$ are computed by employing the
orthogonality of the eigenfunctions, along with the derivative jump
condition
\begin{equation}
\lim_{\varepsilon \to 0} \
{\partial\green \over \partial\chi}\Bigg|_{\chi=\chi_0+\varepsilon}
- {\partial\green \over \partial\chi}\Bigg|_{\chi=\chi_0-\varepsilon}
= - \, {3 \dot N_0 \, \deltapar \, k T_e \, \delta(\tau-\tau_0) \over
\alpha \, \pi r_0^2 \, c \, \epsilon_0^4}
\ ,
\label{eq3.21}
\end{equation}
which is obtained by integrating the transport equation~(\ref{eq3.4})
with respect to $\chi$ in a small region surrounding the injection
energy $\chi_0$. Substituting using the expansion for $\green$ yields
\begin{equation}
\lim_{\varepsilon \to 0} \
\sum_{n=0}^\infty \ C_n \, g_n(\tau) \, \left[h_n'(\chi_0+\varepsilon)
- h_n'(\chi_0-\varepsilon)\right]
= - \, {3 \dot N_0 \, \deltapar \, k T_e \, \delta(\tau-\tau_0) \over
\alpha \, \pi r_0^2 \, c \, \epsilon_0^4}
\ ,
\label{eq3.22}
\end{equation}
where primes denote differentiation with respect to $\chi$. By employing
(\ref{eq3.18}) for $h_n$, we find that
\begin{equation}
\sum_{n=0}^\infty \ C_n \, g_n(\tau) \, \chi_0^{\kappa-4}
\, e^{-\chi_0/2} \mathfrak W(\chi_0)
= - \, {3 \dot N_0 \, \deltapar \, k T_e \, \delta(\tau-\tau_0) \over
\alpha \, \pi r_0^2 \, c \, \epsilon_0^4}
\ ,
\label{eq3.23}
\end{equation}
where the Wronskian, $\mathfrak W$, is defined by
\begin{equation}
\mathfrak W(\chi_0) \equiv M_{\kappa,\mu}(\chi_0) \,
W'_{\kappa,\mu}(\chi_0) - W_{\kappa,\mu}(\chi_0) \,
M'_{\kappa,\mu}(\chi_0)
\ .
\label{eq3.24}
\end{equation}

We can evaluate $\mathfrak W(\chi_0)$ analytically using
equation~(2.4.27) from Slater$^{15}$, which yields
\begin{equation}
\mathfrak W(\chi_0) = - \, {\Gamma(1+2\mu) \over \Gamma(\mu-\kappa+{1 \over 2})}
\ .
\label{eq3.25}
\end{equation}
Using this result to substitute for $\mathfrak W(\chi_0)$ in
(\ref{eq3.23}) and reorganizing the terms, we obtain
\begin{equation}
\sum_{n=0}^\infty \ {\Gamma(1 + 2 \mu) \, C_n \, g_n(\tau) \over
\Gamma(\mu-\kappa+{1 \over 2})}
= {3 \dot N_0 \, \deltapar \, e^{\chi_0/2} \, \delta(\tau-\tau_0) \over
\alpha \, \pi r_0^2 \, c \, \chi_0^\kappa \, (k T_e)^3}
\ ,
\label{eq3.26}
\end{equation}
where $\mu$ is a function of $\lambda_n$ via (\ref{eq3.17}). We can now
calculate the expansion coefficients $C_n$ by utilizing the
orthogonality of the spatial eigenfunctions $g_n$ represented by
(\ref{eq3.15}). Multiplying both sides of (\ref{eq3.26}) by
$e^{3\alpha\tau^2/2} g_m(\tau)$ and integrating with respect to $\tau$
from zero to infinity yields, after some algebra,
\begin{equation}
C_n = {3 \dot N_0 \, \deltapar \, \sqrt{2 w}
\ e^{\chi_0/2} \, e^{3\alpha\tau_0^2/2}
\over \pi r_0^2 \, c (k T_e)^3
\chi_0^\kappa \, \sqrt{\alpha}}
\, {\Gamma(\mu-\kappa+{1 \over 2}) \ n! \,  g_n(\tau_0)
\over \Gamma(1+2\mu) \, \Gamma(n+{1 \over 2})}
\ .
\label{eq3.27}
\end{equation}
The final closed-form solution for the Green's function, obtained by
combining (\ref{eq3.18}), (\ref{eq3.20}), and (\ref{eq3.27}),
is given by
\begin{eqnarray}
\green(\tau_0,\tau,\chi_0,\chi)
&=& {3 \dot N_0 \, \deltapar \, e^{3\alpha\tau_0^2/2} \sqrt{2 w}
\ \chi^{\kappa-4} \, e^{(\chi_0-\chi)/2} \over
\pi r_0^2 \, c (k T_e)^3 \chi_0^\kappa \sqrt{\alpha}}
\sum_{n=0}^\infty \
{\Gamma(\mu-\kappa+{1 \over 2}) \, n! \over \Gamma(1+2\mu) \,
\Gamma(n+{1 \over 2})}
\nonumber
\\
&\times& g_n(\tau_0) \, g_n(\tau)
\, M_{\kappa,\mu}(\chimin)
\, W_{\kappa,\mu}(\chimax)
\ ,
\label{eq3.28}
\end{eqnarray}
where the spatial eigenfunctions $g_n$ are computed using (\ref{eq3.14})
and the parameters $\kappa$ and $\mu$ are given by (\ref{eq3.17}). This
exact, analytical solution for $\green$ provides a very efficient means
for computing the steady-state Green's function resulting from the
continual injection of monochromatic seed photons from a source at an
arbitrary location inside the accretion column. The eigenfunction
expansion converges rapidly and therefore we can generally obtain an
accuracy of at least four significant figures in our calculations of
$\green$ by terminating the series after the first 5-10 terms.

\bigskip
\begin{figure}[t]
\includegraphics[width=100mm]{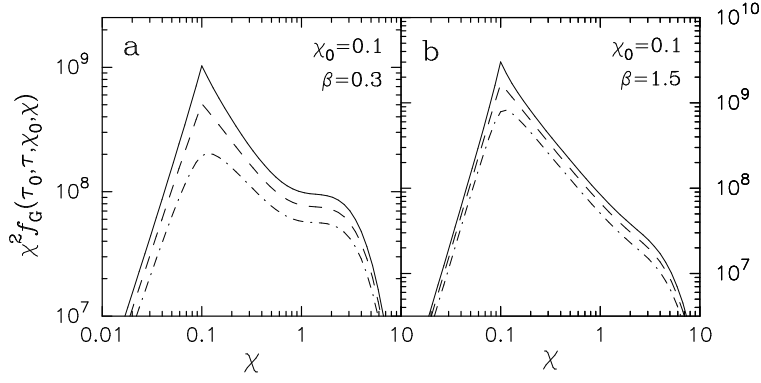}
\caption{Green's function $\chi^2\green$ describing the photon spectrum
inside an x-ray pulsar accretion column [Eq.~(\ref{eq3.28})] plotted
as a function of the dimensionless photon energy $\chi$ for $\dot
N_0=1$, $T_e=10^7\,$K, $r_0=10^4\,$cm, $\alpha=0.1$, $\xi=1.5$,
$\chi_0=0.1$, $\tau_0=0.5$, and (a) $\deltapar=0.3$, (b)
$\deltapar=1.5$. The values of the optical depth are $\tau=0.01$ (solid
line), $\tau=1.0$ (dashed line), and $\tau=1.5$ (dot-dashed line).}
\end{figure}

\bigskip
\begin{figure}[t]
\includegraphics[width=100mm]{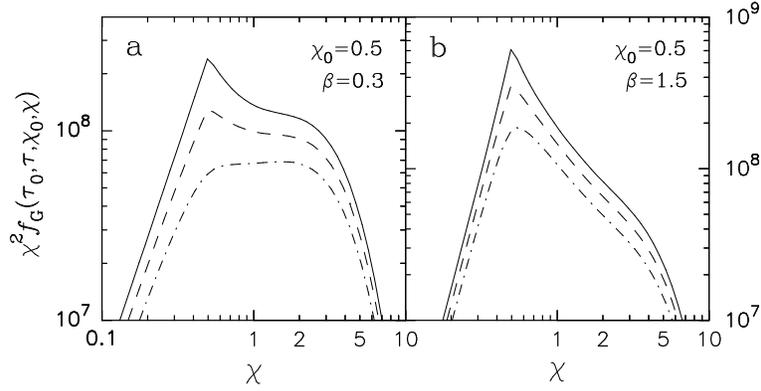}
\caption{Same as Fig.~1, except $\chi_0=0.5$. In this case the energy of
the injected photons is comparable to the thermal energy of the electrons.}
\end{figure}

\bigskip

\section*{\bf 3.5 Numerical Examples}

\bigskip

In this section we illustrate the computational method by examining the
dependence of the Green's function $\green(\tau_0,\tau,\chi_0,\chi)$ on
the optical depth $\tau$ and the dimensionless energy $\chi$. We remind
the reader that the solution for the Green's function represents the
photon spectrum inside the accretion column at the specified position
and energy, resulting from the injection of monochromatic photons with
dimensionless energy $\chi_0$ from a source located at optical depth
$\tau_0$. Analysis of $\green$ therefore allows us to explore the
effects of bulk and thermal Comptonization as the photons diffuse
throughout the plasma, eventually escaping through the walls of the
accretion column. In terms of $\chi$ and $\tau$, the radiation number
density associated with the Green's function is given by
[cf.~Eq.~(\ref{eq2.2})]
\begin{equation}
\ngreen(\tau) = (k T_e)^3 \int_0^\infty \chi^2 \,
\green(\tau_0,\tau,\chi_0,\chi) \, d\chi
\ .
\label{eq3.29}
\end{equation}
It follows that $(k T_e)^3\chi^2\green\,d\chi$ equals the number of
photons per unit volume at optical depth $\tau$ with dimensionless
energy between $\chi$ and $\chi+d\chi$.

In Fig.~1 we use (\ref{eq3.28}) to plot $\chi^2\green$ as a function of
$\chi$ and $\tau$ for the parameter values $\dot N_0=1$, $T_e=10^7\,$K,
$r_0=10^4\,$cm, $\alpha=0.1$, $\xi=1.5$, $\chi_0=0.1$, $\tau_0=0.5$, and
either $\deltapar=0.3$ or $\deltapar=1.5$. We remind the reader that the
value of $\deltapar$ determines the relative importance of bulk and
thermal Comptonization [see Eq.~(\ref{eq3.8})]. Both spectra in Fig.~1
display a peak at the energy of the injected photons, $\chi=0.1$. For
the case with $\deltapar=0.3$, there is also a noticeable Wien hump in
the spectrum around $\chi\sim 2$ due to the effect of thermal
Comptonization.$^{3}$ However, when $\deltapar=1.5$, the Wien feature is
completely hidden by the power-law shape associated with the dominant
bulk Comptonization process. The Green's function is plotted in Fig.~2
for the same parameters used in Fig.~1, except that we now set
$\chi_0=0.5$, so that the energy of the injected photons, $\epsilon_0$,
is comparable to the thermal energy of the electrons, $k T_e$. In this
case, the Wien hump is still visible when $\beta=0.3$, although it has
begun to merge into the peak associated with the injected photons
because of the larger value of $\chi_0$. Note that all of the spectra in
Figs.~1 and 2 display a normalization that decreases with increasing
$\tau$, due to the fact that the downward advection of the electrons
tends to ``trap'' the photons near the bottom of the accretion column,
close to the stellar surface.$^{8}$

Equation~(\ref{eq3.28}) can also be used to calculate the spectrum of
the radiation escaping through the walls of the accretion column, which
comprises the radiation spectrum observed at Earth. Because the x-ray
pulsars are distant cosmic sources, the observed spectra are the result
of emission escaping over the {\it entire} vertical extent of the
accretion column, and therefore we must perform a vertical integration
with respect to $\tau$ in order to compare the model spectra with actual
x-ray observations. This procedure is carried out in reference
[4], where the resulting spectra are compared with the x-ray data
for several astrophysical sources, confirming excellent agreement.

\bigskip

\section*{\bf IV. IDENTITIES INVOLVING THE LAGUERRE POLYNOMIALS}

\bigskip

We are primarily interested here in the mathematical properties of the
solution for the Green's function given by (\ref{eq3.28}). In the
present application, it is possible to compute the photon number
density, $\ngreen$, either by integrating the series expansion for the
Green's function term-by-term, or by solving directly the ordinary
differential equation satisfied by $\ngreen$. By equating the results
obtained for the photon number density using the two methods, we can
derive an interesting new summation identity for the Laguerre
polynomials $L_n^{(-1/2)}(x)$ appearing in the expression for the
spatial eigenfunctions $g_n(\tau)$. We carry out this procedure below.

\bigskip

\section*{\bf 4.1 Integration of the Green's Function}

\bigskip

By using (\ref{eq3.28}) to substitute for $\green$ in (\ref{eq3.29}) and
integrating term-by-term, we find that the associated photon number
density $\ngreen$ can be written as
\begin{equation}
\ngreen(\tau)={3 \dot N_0 \, \deltapar \, e^{3\alpha\tau_0^2/2} \sqrt{2 w}
\ e^{\chi_0/2} \over \pi r_0^2 \, c \, \chi_0^\kappa \sqrt{\alpha}}
\sum_{n=0}^\infty \
{\Gamma(\mu-\kappa+{1 \over 2}) \, n! \over \Gamma(1+2\mu) \,
\Gamma(n+{1 \over 2})} \ g_n(\tau_0) \, g_n(\tau) \, K_{_{\rm G}}(\chi_0)
\ ,
\label{eq4.1}
\end{equation}
where
\begin{equation}
K_{_{\rm G}}(\chi_0) \equiv
W_{\kappa,\mu}(\chi_0) \int_0^{\chi_0}
\chi^{\kappa-2} \, e^{-\chi/2} \, M_{\kappa,\mu}(\chi) \, d\chi
+ M_{\kappa,\mu}(\chi_0) \int_{\chi_0}^\infty
\chi^{\kappa-2} \, e^{-\chi/2} \, W_{\kappa,\mu}(\chi) \, d\chi
\ .
\label{eq4.2}
\end{equation}
The two indefinite integrals on the right-hand side of (\ref{eq4.2}) can be
evaluated using equations~(3.2.6) and (3.2.12) from Slater$^{15}$. After some
algebra, we find that
\begin{equation}
\int_0^{\chi_0} \chi^{\kappa-2} \, e^{-\chi/2} \, M_{\kappa,\mu}(\chi) \, d\chi
= {e^{-\chi_0/2} \, \chi_0^{\kappa-1} \over \kappa+\mu-{1 \over 2}}
\ M_{\kappa-1,\mu}(\chi_0)
\ ,
\label{eq4.3}
\end{equation}
and
\begin{equation}
\int_{\chi_0}^\infty \chi^{\kappa-2} \, e^{-\chi/2} \, W_{\kappa,\mu}(\chi)
\, d\chi = e^{-\chi_0/2} \, \chi_0^{\kappa-1} \ W_{\kappa-1,\mu}(\chi_0)
\ .
\label{eq4.4}
\end{equation}
Combining relations yields
\begin{equation}
K_{_{\rm G}}(\chi_0) = {e^{-\chi_0/2} \, \chi_0^{\kappa-1} \over
\kappa+\mu-{1 \over 2}} \ W_{\kappa,\mu}(\chi_0) \, M_{\kappa-1,\mu}(\chi_0)
+ e^{-\chi_0/2} \, \chi_0^{\kappa-1} \ M_{\kappa,\mu}(\chi_0)
\, W_{\kappa-1,\mu}(\chi_0)
\ ,
\label{eq4.5}
\end{equation}
or, equivalently,
\begin{equation}
K_{_{\rm G}}(\chi_0) = {\chi_0^\kappa \, e^{-\chi_0/2} \over
\kappa+\mu-{1 \over 2}} \left[M'_{\kappa-1,\mu}(\chi_0) \, W_{\kappa-1,\mu}(\chi_0)
- W'_{\kappa-1,\mu}(\chi_0) \, M_{\kappa-1,\mu}(\chi_0)\right]
\ ,
\label{eq4.6}
\end{equation}
where we have used equations~(13.4.32) and (13.4.33) from Abramowitz and
Stegun,$^{13}$ and primes denote differentiation with respect to
$\chi_0$. The Wronskian on the right-hand side of (\ref{eq4.6}) can be
evaluated using (\ref{eq3.25}) to obtain
\begin{equation}
K_{_{\rm G}}(\chi_0) =
{\chi_0^\kappa \, e^{-\chi_0/2} \over \kappa+\mu-{1 \over 2}} \
{\Gamma(1+2\mu) \over \Gamma(\mu-\kappa+{3 \over 2})}
\ ,
\label{eq4.7}
\end{equation}
which can be further simplified by applying (\ref{eq3.17}), yielding
\begin{equation}
K_{_{\rm G}}(\chi_0) = {\chi_0^\kappa \, e^{-\chi_0/2}
\over \deltapar \, (\lambda_n-3)} \ {\Gamma(1+2\mu)
\over \Gamma(\mu-\kappa+{1 \over 2})}
\ ,
\label{eq4.8}
\end{equation}
where the eigenvalues $\lambda_n$ are given by (\ref{eq3.13}). By using
(\ref{eq4.8}) to substitute for $K_{_{\rm G}}(\chi_0)$ in (\ref{eq4.1})
and simplifying, we find that the photon number density can be evaluated
using the expansion
\begin{equation}
\ngreen(\tau)={3 \dot N_0 \, e^{3\alpha\tau_0^2/2} \sqrt{2 w}
\over \pi r_0^2 \, c \, \sqrt{\alpha}}
\sum_{n=0}^\infty \
{n! \, g_n(\tau_0) \, g_n(\tau) \over \Gamma(n+{1 \over 2}) \, (\lambda_n-3)}
\ .
\label{eq4.9}
\end{equation}
This expression allows the computation of the photon number density at
any optical depth $\tau$ inside the accretion column for given values of
the photon injection rate $\dot N_0$ and the source location $\tau_0$.

\bigskip

\section*{\bf 4.2 Solution of the Differential Equation}

\bigskip

The availability of the transport equation~(\ref{eq3.4}) provides us
with an alternative means for computing the photon number density,
$\ngreen$, by directly solving the differential equation
\begin{equation}
{1 \over 3}\,{d^2\ngreen\over d\tau^2}
+ \alpha \, \tau \, {d\ngreen \over d\tau}
+ \alpha \, \ngreen
- \xi^2 \alpha^2\,\tau^2 \, \ngreen
= 
- {\dot N_0 \, \delta(\tau-\tau_0)
\over \pi r_0^2 c}
\ ,
\label{eq4.10}
\end{equation}
which is derived by operating on (\ref{eq3.4}) with $(k T_e)^3
\int_0^\infty \chi^2\,d\chi$. The homogeneous version of (\ref{eq4.10})
obtained when $\tau\ne\tau_0$ admits the fundamental solutions
\begin{equation}
\ngreen(\tau) =
\begin{cases}
A \, e^{-\alpha(3+w)\tau^2/4} \, U\big(a,{1 \over 2},{\alpha w \tau^2
\over 2}\big) \ , & \tau > \tau_0 \ , \cr
B \, e^{-\alpha(3+w)\tau^2/4} \, M\big(a,{1 \over 2},{\alpha w \tau^2
\over 2}\big) \ , & \tau < \tau_0 \ , \cr
\end{cases}
\label{eq4.11}
\end{equation}
where $A$ and $B$ are constants and
\begin{equation}
a \equiv {w - 3 \over 4 w} \ ,
\ \ \ \ \
w \equiv \left(9 + 12 \, \xi^2\right)^{1/2}
\ .
\label{eq4.12}
\end{equation}
The solutions in (\ref{eq4.11}) are consistent with the physical
boundary conditions at large and small $\tau$ discussed in \S~3.1 [see
Eq.~(\ref{eq3.9})]. The constants $A$ and $B$ are determined by
requiring that $\ngreen(\tau)$ be continuous at the source location
$\tau_0$, and that it satisfy the derivative jump condition
\begin{equation}
\lim_{\varepsilon \to 0} \
{d\ngreen \over d\tau}\Bigg|_{\tau=\tau_0+\varepsilon}
- {d\ngreen \over d\tau}\Bigg|_{\tau=\tau_0-\varepsilon}
= - \, {3 \dot N_0 \over \pi r_0^2 \, c}
\ ,
\label{eq4.13}
\end{equation}
obtained by integrating (\ref{eq4.10}) with respect to $\tau$ in a small
region around $\tau_0$. Combining (\ref{eq4.11}) and (\ref{eq4.13}) with
the continuity condition and the expression for the Wronskian given by
(\ref{eq3.12}), we find after some algebra that
\begin{eqnarray}
A &=& {3 \dot N_0 \Gamma(a) \over \pi r_0^2 c
\sqrt{2 \pi \alpha w}} \ e^{\alpha(3-w)\tau_0^2/4}
M\left(a,{1 \over 2},{\alpha w \tau_0^2 \over 2}\right) \ ,
\\
B &=& {3 \dot N_0 \Gamma(a) \over \pi r_0^2 c
\sqrt{2 \pi \alpha w}} \ e^{\alpha(3-w)\tau_0^2/4} \,
U\left(a,{1 \over 2},{\alpha w \tau_0^2 \over 2}\right) \ ,
\label{eq4.14}
\end{eqnarray}
and consequently the global solution for the photon number density
is given by
\begin{equation}
\ngreen(\tau) = {3 \dot N_0 \Gamma(a) \over \pi r_0^2 c
\sqrt{2 \pi \alpha w}} \ e^{\alpha(3-w)\tau_0^2/4}
e^{-\alpha(3+w)\tau^2/4} \, M\left(a,{1 \over 2},
{\alpha w \taumin^2 \over 2}\right) \, U\left(a,{1 \over 2},
{\alpha w \taumax^2 \over 2}\right)
\,
\label{eq4.15}
\end{equation}
where
\begin{equation}
\taumin \equiv \min(\tau,\tau_0) \ ,
\ \ \ \ \
\taumax \equiv \max(\tau,\tau_0)
\ .
\label{eq4.16}
\end{equation}

\bigskip

\section*{\bf 4.3 Summation Identity for the Laguerre Polynomials}

\bigskip

Equations~(\ref{eq4.9}) and (\ref{eq4.15}) provide two independent
expressions that can each be used to calculate the photon number density
$\ngreen$. This fact allows us to derive an interesting new summation
identity involving the Laguerre polynomials $L_n^{(-1/2)}(x)$ appearing
in equation~(\ref{eq3.14}) for the spatial eigenfunctions $g_n$. By
setting (\ref{eq4.9}) and (\ref{eq4.15}) equal and substituting for
$\lambda_n$ using (\ref{eq3.13}), we obtain after some simplification
\begin{equation}
\sum_{n=0}^\infty {n!
\, L_n^{(-1/2)}({\alpha w \tau_0^2 \over 2})
\, L_n^{(-1/2)}({\alpha w \tau^2 \over 2})
\over (n+a) \, \Gamma(n+{1 \over 2})}
= {\Gamma(a) \over \sqrt{\pi}}
\ M\left(a,{1 \over 2},
{\alpha w \taumin^2 \over 2}\right) U\left(a,{1 \over 2},
{\alpha w \taumax^2 \over 2}\right)
\ .
\label{eq4.17}
\end{equation}
Setting $x=\alpha w \tau^2/2$ and $x_0=\alpha w \tau_0^2/2$ yields the
more compact form
\begin{equation}
\sum_{n=0}^\infty {n! \, L_n^{(-1/2)}(x_0) \, L_n^{(-1/2)}(x)
\over (n+a) \, \Gamma(n+{1 \over 2})}
= {\Gamma(a) \over \sqrt{\pi}}
\ M\left(a,{1 \over 2},\xmin\right) U\left(a,{1 \over 2},\xmax\right)
\ ,
\label{eq4.18}
\end{equation}
where
\begin{equation}
\xmin \equiv \min(x,x_0) \ ,
\ \ \ \ \
\xmax \equiv \max(x,x_0)
\ .
\label{eq4.19}
\end{equation}
Equation~(\ref{eq4.18}) has not appeared in the previous literature,
although Exton,$^{17}$ Srivastava,$^{18}$ and Manocha$^{19}$ have
derived several related identities.

\section*{\bf V. CONCLUSION}

In this article we have applied the methods of classical analysis to
derive the closed-form solution for the Green's function, $\green$,
given by (\ref{eq3.28}), describing the bulk and thermal Comptonization
of monochromatic seed photons scattered by hot, infalling electrons
inside an x-ray pulsar accretion column. X-ray pulsars are rotating
neutron stars, which are the most dense and compact solid objects known
to exist in the universe. These enigmatic sources are characterized by
super strong magnetic, gravitational, and radiation fields, making them
the most extreme physical ``laboratories'' in the universe. It is
therefore of great theoretical interest to obtain the best possible
understanding of the spectral formation process in x-ray pulsars,
because the observed radiation provides us with the only available
window into their physical nature. As demonstrated in Figs.~1 and 2, the
Green's function is characterized by a power-law shape at moderate
photon energies, with an exponential turnover at higher energies, in
agreement with the spectra of many sources.$^{4}$ We conclude that bulk
and thermal Comptonization in the shocked pulsar accretion column
provides a natural explanation for the typical high-energy spectra
produced by x-ray pulsars.

Our primary interest here is in exploring the mathematical properties of
the Green's function. In particular, we established that by performing
two independent calculations of the photon number density $\ngreen$
associated with $\green$, one can obtain an interesting new summation
identity involving the Laguerre polynomials, $L_n^{(-1/2)}(x)$,
expressed by (\ref{eq4.18}). The derivation is based on the simultaneous
calculation of $\ngreen$ using either term-by-term integration of the
Green's function expansion (\ref{eq3.28}), which yields (\ref{eq4.9}),
or instead via the direct solution of the differential equation for
$\ngreen$ given by (\ref{eq4.10}), which leads to (\ref{eq4.15}). The
new identity for the Laguerre polynomials obtained here is related to
various similar expressions presented by Chatterjea,$^{16}$
Exton,$^{17}$ Srivastava,$^{18}$ Manocha,$^{19}$ Varma,$^{20}$ and
Gradshteyn and Ryzhik,$^{14}$ although our results are distinct from
theirs. We expect that our new expression may be of potential benefit in
various problems of mathematical physics.

The author would like to gratefully acknowledge the helpful comments
provided by the anonymous referee.

\section*{REFERENCES}

\smallskip\noindent
$^{1}$P. A. Becker, J. Math. Phys. {\bf 46}, 53511 (2005).

\smallskip\noindent
$^{2}$P. A. Becker and M. T. Wolff, Astrophys. J. {\bf 630}, 465 (2005).

\smallskip\noindent
$^{3}$G. B. Rybicki and A. P. Lightman, {\it Radiative Processes in Astrophysics}
(Wiley, New York, 1979).

\smallskip\noindent
$^{4}$M. T. Wolff, P. A. Becker, and K. D. Wolfram, in {\it The Multicoloured
Landscape of Compact Objects and their Explosive Progenitors},
Cefalu, Sicily, June 2006, ed. L. Burderi et al. (New York: AIP).

\smallskip\noindent
$^{5}$R. D. Blandford and D. G. Payne, Mon. Not. R. Astron. Soc. {\bf 194},
1033 (1981).

\smallskip\noindent
$^{6}$P. A. Becker, Mon. Not. R. Astr. Soc. {\bf 343}, 215--240 (2003).

\smallskip\noindent
$^{7}$Yu. E. Lyubarskii and R. A. Sunyaev, Soviet Astr. Lett. {\bf 8}, 330 (1982).

\smallskip\noindent
$^{8}$P. A. Becker, Astrophys. J. {\bf 498}, 790 (1998).

\smallskip\noindent
$^{9}$R. D. Blandford and D. G. Payne, Mon. Not. R. Astron. Soc.
{\bf 194}, 1041 (1981).

\smallskip\noindent
$^{10}$D. G. Payne and R. D. Blandford, Mon. Not. R. Astron. Soc.
{\bf 196}, 781 (1981).

\smallskip\noindent
$^{11}$M. Colpi, Astrophys. J. {\bf 326}, 223 (1988).

\smallskip\noindent
$^{12}$P. Schneider and J. G. Kirk, Astrophys. J. Lett. {\bf 323}, L87 (1987).

\smallskip\noindent
$^{13}$M. Abramowitz and I. A. Stegun, {\it Handbook of Mathematical
Functions} (Dover, New York, 1970).

\smallskip\noindent
$^{14}$I. S. Gradshteyn and I. M. Ryzhik, {\it Table of Integrals,
Series, and Products} (Academic Press, London, 1980).

\smallskip\noindent
$^{15}$L. J. Slater, {\it Confluent Hypergeometric Functions} (Cambridge
Univ. Press, Cambridge, England, 1960).

\smallskip\noindent
$^{16}$S. K. Chatterjea, Ann. Scuola Norm. Sup. Pisa {\bf 20}, 739 (1966).

\smallskip\noindent
$^{17}$H. Exton, J. Comp. Appl. Math. {\bf 100}, 225 (1998).

\smallskip\noindent
$^{18}$H. M. Srivastava, Panamer. Math. J. {\bf 1}, 69 (1991).

\smallskip\noindent
$^{19}$B. L. Manocha, Bull. Math. Soc. Sci. Math. R. S. Roumanie
{\bf 11}, 85 (1967).

\smallskip\noindent
$^{20}$V. K. Varma, J. Indian Math. Soc. {\bf 32}, 1 (1968).

\end{document}